\documentclass[12pt]{amsart}

\usepackage{amsthm, amssymb, amsmath}

\theoremstyle{plain}
\numberwithin{equation}{section}
\newtheorem{theorem}[equation]{Theorem}
\newtheorem{prop}[equation]{Proposition}
\newtheorem{lemma}[equation]{Lemma}
\newtheorem{claim}[equation]{Claim}
\newtheorem{observation}[equation]{Observation}

\theoremstyle{definition}

\title{Function Theory on the Neil Parabola}

\author{Greg Knese} \address{Department of Mathematics \\ Washington
University in St. Louis \\ St. Louis, Missouri 63130}
\email{geknese(AT)math.wustl.edu}
\urladdr{http://math.wustl.edu/\textasciitilde geknese}

\date{August 18, 2006}

\begin{document}

\begin{abstract}
  We give a formula for the Carath\'eodory distance on the Neil
  Parabola $\{(z,w) \in \mathbb{D}^2: z^2=w^3 \}$ restricted to the
  bidisk, making it the first variety with a singularity to have its
  Carath\'eodory pseudo-distance explicitly computed.  This addresses
  a recent question of Jarnicki and Pflug. In addition, we relate this
  problem to a mixed Carath\'eodory-Pick interpolation problem for
  which known interpolation theorems do not apply.  Finally, we prove
  a bounded holomorphic function extension result from the Neil
  parabola to the bidisk.
\end{abstract}

\maketitle

\begin{section}{Introduction}
  Distances on a complex space $X$ which are invariant under
  biholomorphic maps have played an important role in the geometric
  approach to complex analysis.  One of the oldest such distances is
  the the Carath\'eodory pseudodistance $c_X$ (``pseudo'' because the
  distance between two points can be zero).  It was introduced by C.\!
  Carath\'eodory in 1926 and is extremely simple to define.  The
  distance between two points $x$ and $y$ is defined to be the largest
  distance (using the Poincar\'e hyperbolic distance) that can occur
  between $f(x)$ and $f(y)$ under a holomorphic map $f$ from $X$ to
  the unit disk $\mathbb{D} \subset \mathbb{C}$.  The Kobayashi
  pseudodistance $k_X$, introduced by S.\! Kobayashi in 1967, is
  defined in the opposite direction: the ``distance'' between two
  points $x$ and $y$ is now the infimum of the (hyperbolic) distance
  that can occur between two points $a, b \in \mathbb{D}$ for which
  there is a holomorphic map $f$ from the disk to $X$ mapping $a$ to
  $x$ and $b$ to $y$.  (Actually, there is a small technicality
  here---see section 4 for the true definition). A consequence of the
  Schwarz-Pick lemma on the disk (which says holomorphic self-maps of
  the disk are distance decreasing in the hyperbolic distance) is the
  fact that $c_X \leq k_X$.
  
  For the purposes of motivating the present article, let us indulge
  in a short tangent.  An interesting question, because of its
  geometric implications (including the existence of one dimensional
  analytic retracts), is for which complex spaces do we have $c_X =
  k_X$?  The most important contribution to this question is by L.\!
  Lempert \cite{LEM}.  Lempert's theorem proves the Carath\'eodory and
  Kobayashi distances agree on a convex domain.  This theorem came as
  a surprise for a couple of reasons: first, convexity is not a
  biholomorphic invariant, and second, which is our main point here,
  \emph{there were not many explicit examples available at the
  time}\footnote{The plot thickens on this problem: there is a domain,
  namely the symmetrized bidisc, in $\mathbb{C}^2$ for which the two
  distances agree, yet this domain is not biholomorphically equivalent
  to a convex domain.  See \cite{J-P} for a summary of these
  results.}.  While we cannot remedy the problem of a lack of examples
  in the past, we can attempt to add to the current selection of
  explicit examples.  Many theorems about invariant metrics can be
  proved in the generality of complex spaces (see \cite{KOB} for
  instance) yet curiously there do not seem to be any nontrivial,
  explicit examples of the Carath\'eodory distance for a complex space
  \emph{with a singularity}.  Perhaps the simplest complex space with
  a singularity is the variety contained in the bidisk given by
\[
N = \{(z,w) \in \mathbb{D}^2: z^2=w^3\}
\]
Following \cite{J-P}, we shall call this the \emph{Neil
  parabola}\footnote{The real curve $y^2=x^3$ and its variations are
  referred to as Neil's semi-cubical parabola. Named after William
  Neil (sometimes spelled Neile), a student of John Wallis, it was the
  first algebraic curve to have its arc length computed via
  proto-calculus techniques \cite{WAL}.}.  In their recent follow-up
  \cite{J-P} to their monograph \cite{J-P93}, M.\! Jarnicki and P.\!
  Pflug pose the following problem:

\emph{Find an effective formula for the
  Carath\'eodory distance on the Neil parabola $N$.}

  In this paper, we give an answer to this problem (see theorem
  \ref{MainThm}).  In addition, we compute the infinitesimal
  Carath\'eodory pseudodistance for the Neil parabola (see theorem
  \ref{thm2}).  As applications, we prove a mixed Carath\'eodory-Pick
  interpolation result for which known interpolation theorems do not
  apply (see theorem \ref{prop:interp}) and we prove a result on
  extending bounded holomorphic functions on the Neil parabola to the
  entire bidisk (see theorem \ref{thm:extend}).

The general layout of the rest of the paper is as follows.  Motivation
  and background for the two previously mentioned applications are
  presented in the following two subsections. In section \ref{section:def},
  precise statements of definitions and results are given along with a
  subsection on preliminary facts about complex analysis on the Neil
  parabola.  The rest of the paper is devoted to proofs. (The
  locations of specific proofs are given near the corresponding
  theorem statements in section \ref{section:def}.)

%%%%%%%%%%%%%%%%%%%%%%%%%%%%%%%%%%%%%%%%%%%%%%%%%%%%%%%%%%%%%%%%%%%%

\begin{subsection}{A mixed Carath\'eodory-Pick problem}\label{section:interp}
 Given $n$ points in the unit disk $z_i$ and $n$ target values $w_i$
  also in the unit disk, the well-known theorem of G.\! Pick
  \cite{PICK} says exactly when there exists a holomorphic
  $F:\mathbb{D} \to \mathbb{D}$ satisfying $F(z_i)=w_i$ (this problem
  was studied independently by Nevanlinna \cite{NEV}).  In fact, the
  Schwarz-Pick lemma is just the version of this for two points: $z_1,
  z_2$ can be interpolated to $w_1, w_2$ if and only if
\[
\left| \frac{w_1-w_2}{1-\bar{w_1} w_2} \right| \leq \left|
  \frac{z_1-z_2}{1-\bar{z_1} z_2} \right|
\]
Similarly, given $n$ complex numbers $a_0, a_1, \dots, a_{n-1}$ a
well-known theorem of Carath\'eodory and Fej\'er \cite{CF} says when
there exists a holomorphic function $F:\mathbb{D}\to
\overline{\mathbb{D}}$ with $a_0, a_1, \dots, a_{n-1}$ as the first
$n$ Taylor coefficients of $F$.\footnote{Using $\overline{\mathbb{D}}$
instead of $\mathbb{D}$ is just a trick used to include the constant
unimodular valued functions, because we are really talking about
functions in the closed unit ball of $H^{\infty}(\mathbb{D})$. The
same idea applies later on to $\mathcal{O}(\mathbb{D},
\overline{\mathbb{D}})$ (of course this notation has not yet been
introduced).}  For $n=2$, this is given again by the (infinitesimal)
Schwarz-Pick lemma: $a_0$ and $a_1$ can be the first two Taylor
coefficients exactly when
\[
|a_0|^2 + |a_1| \leq 1
\]
The first kind of interpolation problem above is called
Nevanlinna-Pick interpolation and the second is called
Carath\'eodory-Fej\'er interpolation.  More modern proofs of these
theorems, using ideas from operator theory like the commutant lifting
theorem of Sz.-Nagy and Foia\c{s} and reproducing kernel Hilbert spaces
(see \cite{FF} and \cite{AM}), make it possible to study so-called
\emph{mixed Carath\'eodory-Pick problems} wherein the idea is to
specify several Taylor coefficients at several points in the disk and
determine whether there exists a holomorphic function from the disk to
the disk with those properties.  However, a restriction imposed in all
of the usual mixed Carath\'eodory-Pick problems is that the Taylor
coefficients must be specified sequentially (i.e.\! one cannot ask to
specify the first and third Taylor coefficients at a point without
specifying the second as well).  For example, these problems do not
address an interpolation problem of the following form: given $z_1,
z_2, z_3, w_1, w_2 \in \mathbb{D}$, when is there a holomorphic
function $F:\mathbb{D} \to \mathbb{D}$ satisfying the following?

\begin{align}
F(z_1) & = w_1 \nonumber\\
F(z_2) & = w_2 \label{mixedcp}\\
F'(z_3) & = 0 \nonumber
\end{align}

In fact, as we shall see, solving the problem (\ref{mixedcp}) amounts
to computing the Carath\'eodory distance for the Neil parabola.  See
theorem \ref{prop:interp} for the exact statement of our result.

\end{subsection} % end of cara-pick problem subsection

%%%%%%%%%%%%%%%%%%%%%%%%%%%%%%%%%%%%%%%%%%%%%%%%%%%%%%%%%%%%%%%%%%%%%

\begin{subsection}{Extension of bounded holomorphic functions on the
    Neil parabola}
\label{section:extend}

  The following result is a special case of the work of H. Cartan on
  Stein Varieties (see \cite{GUN3} page 99). (In fact, we are stating
  it in almost as little generality as possible.)
\begin{theorem}[Cartan] Every holomorphic function on a subvariety $V$
  of $\mathbb{D}^2$ is the restriction of a holomorphic function on
  all of $\mathbb{D}^2$.
\end{theorem}

A vast improvement on this theorem (again stated in simple terms) was
given by P.L.\! Polyakov and G.M.\! Khenkin \cite{P-H}.  They proved
using the methods of integral formulas that any subvariety $V$ of
$\mathbb{D}^2$ satisfying a certain transversality condition has the
property that any bounded holomorphic function on $V$ can be extended
to a bounded holomorphic function on all of $\mathbb{D}^2$.  In fact,
there is a bounded linear operator $T: H^{\infty}(V) \to H^{\infty}
(\mathbb{D}^2)$ with $Tf\mid_V = f$; in other words, there is some
constant $C$ such that for any $f \in H^{\infty} (V)$
\begin{equation}
||Tf||_\infty \leq C ||f||_\infty
\label{operator}
\end{equation}
The previously mentioned ``transversality condition'' applies to the
Neil parabola, and therefore any bounded holomorphic function on $N$
can be extended to a bounded holomorphic function on the bidisk.

Related to these ideas is a paper of J.\! Agler and J.E.\! M{\raise
  .45ex\hbox{c}}Carthy \cite{AM-03}, which gives a description of
varieties in the bidisk with the property that bounded holomorphic
functions can be extended to the bidisk without increasing their
$H^{\infty}$ norm.  The Neil parabola is not such a variety as their
results show.  This can be seen relatively easily from the fact that
the Carath\'eodory pseudodistance on the Neil parabola is not the
restriction of the Carath\'eodory pseudodistance on the bidisk.
Meaning, there is some holomorphic function from $N$ to $\mathbb{D}$
which separates two points of $N$ farther than a function from the
bidisk to the disk could.  Hence, such a function could not be
extended to the bidisk without increasing its norm.

This suggests that extremal functions on the Neil parabola for the
Carath\'eodory pseudodistance might be good candidates for functions
which extend ``badly'' to the bidisk.  Indeed, this allows us to give
a lower bound of $5/4$ on the constant $C$ in (\ref{operator}) for the
Neil Parabola.  In addition to this we present a simple proof using
Agler's Nevanlinna-Pick interpolation theorem for the bidisk that any
bounded holomorphic function on the Neil parabola can be extended to
a bounded holomorphic function on the bidisk with norm increasing by
at most a factor of $\sqrt{2}$ if the function vanishes at the origin
and by a factor of $2\sqrt{2}+1$ otherwise.  This does not exactly
reprove Polyakov and Khenkin's result in our context, since we are not
claiming the extension can be given by a linear operator.
Nevertheless, it is certainly relevant to their result, is much easier
to prove, and provides an explicit bound (see theorem
\ref{thm:extend}).

\end{subsection} % end of extension subsection

\end{section} % end of introduction
%%%%%%%%%%%%%%%%%%%%%%%%%%%%%%%%%%%%%%%%%%%%%%%%%%%%%%%%%%%%%%%%%%%%%%%%

\begin{section}{Definitions and statements of results}\label{section:def}
 
  Let us define several important notions for this paper.  We shall
  use $\mathcal{O}(X,Y)$ to denote the set of holomorphic maps from
  $X$ to $Y$ and $\mathcal{O}(X)$ to denote the set of holomorphic
  functions from $X$ to $\mathbb{C}$, where $X$ and $Y$ are complex
  spaces possibly containing singularities (this holds for $X$ below
  as well).

\begin{itemize}
\item Frequent use will be made of the family of holomorphic
  automorphisms $\phi_\alpha$ of the unit disk $\mathbb{D} \subset
  \mathbb{C}$ given by
\begin{equation}
\phi_\alpha (z) = \frac{\alpha-z}{1-\bar{\alpha}z}
\label{phi}
\end{equation}
where $\alpha \in \mathbb{D}$.  Note that $\phi_\alpha$ is its own
inverse function.  Sometimes we allow $\alpha$ to be in $\partial
\mathbb{D}$, but keep in mind that the resulting $\phi_\alpha$ is no
longer an automorphism of the disk and is instead the constant
function $\bar{\alpha}$.

\item The \emph{pseudo-hyperbolic distance} on $\mathbb{D}$ is defined
to be
\[
m(a,b)= \left| \frac{a-b}{1-\bar{a}b} \right|
\]
The \emph{Poincar\'e distance} on $\mathbb{D}$ is given by $\rho=\tanh^{-1}
m$.

\item The \emph{Poincar\'e metric} on the disk, which we shall also
  denote by $\rho$, is defined to be
\[
\rho(z; v) = \frac{|v|}{1-|z|^2}
\]
for $z\in \mathbb{D}$ and $v \in \mathbb{C}$.

\item The \emph{Carath\'eodory pseudodistance} on $X$
  is denoted by $c_X$ and is defined by
\[
c_{X}(x,y) := \sup\{\rho(f(x),f(y)): f \in \mathcal{O}(X,\mathbb{D})\}
\]
If we replace $\rho$ above with $m$, we get what Jarnicki and Pflug
call the \emph{M\"obius pseudodistance}:
\[
c^{*}_X (x,y) := \sup\{m(f(x),f(y)): f \in
\mathcal{O}(X,\mathbb{D})\}
\]

Due to the simple formula for $m$ and the relation $c_X= \tanh^{-1}
c^{*}_X$, the M\"obius pseudodistance is more computationally useful
for our purposes, and therefore will be used exclusively in all
proofs.

\item The \emph{Carath\'eodory pseudometric} $C_X$ is defined to be
\[
C_X (x; v) = \sup\{ \rho(f(x); df_x (v)) : f \in \mathcal{O}(X,
\mathbb{D})\}
\]
for $x \in X$ and $v \in T_x X$, the tangent space of $X$ at $x$.
The Carath\'eodory pseudometric will often be referred to as the
\emph{infinitesimal Carath\'eodory pseudodistance}.

\item Finally, the \emph{Lempert function} for $X$ is
  denoted $\tilde{k}_X$ and is defined by
\[
\tilde{k}_X (x,y) = \inf\{\rho(a,b): \exists f\in
  \mathcal{O}(\mathbb{D}, X) \text{ with } f(a)=x, f(b)=y\}
\]
where $\tilde{k}_X$ is defined to equal $\infty$ if the above set over
which the infimum is taken is empty.  The \emph{Kobayashi
  pseudodistance} $k_X$ is then defined to be largest pseudodistance
bounded by $\tilde{k}_X$.
\end{itemize}

For more information on and examples of these definitions we refer the
reader to \cite{IK}, \cite{J-P93}, \cite{J-P}, and \cite{KOB}.

Recall from the introduction that the Neil parabola is the set
\[
N = \{ (z,w)\in \mathbb{D}^2: z^2=w^3 \}
\]
The set $N$ is a one-dimensional connected analytic variety in
$\mathbb{D}^2$ with a singularity at $(0,0)$.  Furthermore, $N$ has a
bijective holomorphic parametrization $p: \mathbb{D} \to N$ given by
\begin{equation}
p(\lambda) := (\lambda^3, \lambda^2)
\label{param}
\end{equation}
The function $q:=p^{-1}$ is continuous on $N$, holomorphic on
$N\setminus \{(0,0)\}$, and can be given by $q(z,w)=z/w$ when
$(z,w)\ne (0,0)$ (and $q(0,0)=0$).  For the benefit of those readers
unfamiliar with holomorphic functions on a variety with a singularity,
we include a discussion of these ideas in the concrete context of the
Neil parabola in section \ref{section:discussion}. It is known that
the Kobayashi pseudodistance $k_{N}$ and the Lempert function
$\tilde{k}_{N}$ for $N$ are as simple as possible (see \cite{J-P}):
\[
k_N ((a,b),(z,w)) = \tilde{k}_{N} ((a,b),(z,w))= \rho(q(a,b),q(z,w))
\]
On the other hand (and to reiterate our goal in
this paper), in \cite{J-P} the authors lament that despite the
simplicity of $N$ an effective formula for the Carath\'eodory
pseudodistance $c_N$ is not known.  We propose the following as an
effective formula for $c_N$.

\begin{theorem}[Carath\'eodory pseudodistance formula]
  Given nonzero $\lambda, \delta \in \mathbb{D}$, let
\[
\alpha_0=\frac{1}{2}\left(\frac{1}{\bar{\lambda}} +
  \lambda+\frac{1}{\bar{\delta}}+\delta \right)
\]
then
\[
c_N(p(\lambda), p(\delta)) = \begin{cases} \rho(\lambda^2, \delta^2)
  &\text{
    if } |\alpha_0| \geq 1 \\
  \rho(\lambda^2\phi_{\alpha_0} (\lambda), \delta^2\phi_{\alpha_0}
  (\delta)) &\text{ if } |\alpha_0| < 1
    \end{cases}
\]
Also, $c_N (p(0),p(\lambda)) = \rho(0,\lambda^2)=\tanh^{-1} |\lambda|^2$.
\label{MainThm}
\end{theorem}

In particular, it should be noted that if $\lambda$ and $\delta$ have
an acute angle between them (i.e.\! ${\rm Re} (\lambda \bar{\delta}) >
0$), then $|\alpha_0|>1$, and the first formula above gives the
distance between $p(\lambda)$ and $p(\delta)$.  Also, the theorem
shows $k_N \ne c_N$ as one might suspect.

In section \ref{section:reduction} we shall reduce the above problem
to a maximization problem on the closed unit disk, and in section
\ref{section:MainThmProof} we solve the maximization problem to yield
theorem \ref{MainThm}.  In addition, a slightly nicer form of the
above formula will be presented as proposition \ref{prop:nicer}.

As will be explained in subsection \ref{section:discussion}, the
tangent spaces of $N$ can be identified with subspaces of the tangent
spaces of $\mathbb{D}^2$.  In particular, for $x=(a,b) \ne (0,0)$,
$T_x N$ is simply the span of the vector $(3a,2b)$, while the tangent
space at the origin of $N$ is two dimensional and therefore equal to
all of $\mathbb{C}^2=T_{(0,0)} \mathbb{D}^2$.  We can now present our
formula for the Carath\'eodory pseudometric of $N$ (this is proved in
section \ref{section:infproof}).

\begin{theorem}[Carath\'eodory pseudometric formula]
  For $v=(v_1, v_2) \in \mathbb{C}^2$, we have
\begin{equation}
C_N ((0,0); v) = \begin{cases} |v_2| & \text{ if } |v_2| \geq 2 |v_1| \\
                               \frac{4|v_1|^2+|v_2|^2}{4|v_1|} &
                               \text{ if } |v_2| < 2 |v_1| \end{cases}
\label{metric:origin}
\end{equation}
and for $(a,b) \in N$ nonzero and $z \in \mathbb{C}$ we have
\begin{equation}
C_N ((a,b); z(3a,2b)) = \frac{2|b|}{1-|b|^2} |z|
\label{metric:general}
\end{equation}

\label{thm2}
\end{theorem}

As mentioned in subsection \ref{section:interp}, a direct consequence
of the preceding formulas is the following atypical mixed
Carath\'eodory-Pick interpolation result (see section
\ref{section:proveinterp} for the proof).  

\begin{theorem}[Mixed interpolation problem]
  Given distinct $z_1, z_2, z_3 \in \mathbb{D}$ and $w_1, w_2
  \in \mathbb{D}$, there exists $F \in \mathcal{O}(\mathbb{D},
  \mathbb{D})$ with
\[
\begin{aligned}
F(z_i) & = w_i \text{ for } i=1,2\\
F'(z_3) & = 0
\end{aligned}
\]
if and only if
\begin{equation}
\rho (w_1, w_2) \leq c_N (p(\phi_{z_3}(z_1)), p(\phi_{z_3}(z_2)))
\label{mixedcp-ineq}
\end{equation}
Moreover, if the problem is extremal (i.e. if there is equality in
(\ref{mixedcp-ineq})), then the solution is unique and is a Blaschke
product of order two or three.
\label{prop:interp}
\end{theorem}

The significance of the theorem (which one could say in the way it is
stated now practically follows from definitions) is of course that
$c_N$ is directly computable by theorem \ref{MainThm}. (So, inequality
\eqref{mixedcp-ineq} is easy to check.)

Finally, in section \ref{section:proveextend} we prove the following
result on extending bounded holomorphic functions from the Neil
parabola to the bidisk.

\begin{theorem}[Bounded analytic extension]  
  For any $f \in \mathcal{O}(N, \mathbb{D})$ with $f(0,0)=0$, there
  exists an extension of $f$ to a function in
  $\mathcal{O}(\mathbb{D}^2, \sqrt{2}\mathbb{D})$.  If $f(0,0) \ne 0$,
  then $f$ can be extended to $\mathcal{O}(\mathbb{D}^2,
  (2\sqrt2+1)\mathbb{D})$.  In addition, there exists a function in
  $\mathcal{O}(N, \mathbb{D})$ which cannot be extended to a function
  in $\mathcal{O}(\mathbb{D}^2, r\mathbb{D})$ for $r < 5/4$.
\label{thm:extend}
\end{theorem}

Here $r\mathbb{D}$ just refers to the disk of radius $r$.

\begin{subsection}{Complex analysis on the Neil parabola}\label{section:discussion}
  In this subsection we discuss how to do complex analysis on a
  variety with a singularity in the concrete setting of the Neil
  parabola.  This is adapted from \cite{J-P} and \cite{GUN} (see pages
  18-20 and the chapter on tangent spaces) and nothing in this section
  is by any means new.  The most important facts of this subsection
  are summarized in the two ``observations'' \ref{observation1} and
  \ref{observation2}.

A function $f$ on $N$ is defined to be holomorphic if at each point
  $x\in N$, there is a holomorphic function $F$ on a neighborhood $U$
  of $x$ in the bidisk which agrees with $f$ on $U\cap
  N$. Fortunately, we can give a more concrete description of the set
  of holomorphic functions on $N$.  Given $f \in \mathcal{O}(N)$, the
  function $h:=f \circ p$ (recall $p$ from \eqref{param}) is an element
  of $\mathcal{O}(\mathbb{D})$ satisfying $h'(0)=0$.  The reason for
  this is given an extension, $F$, of $f$ holomorphic on a
  neighborhood of $(0,0)$ in $\mathbb{D}^2$, $h=F \circ p$ is
  holomorphic on a neighborhood of $0$ in $\mathbb{D}$.  Hence, the
  derivative $h'(\lambda)=dF_{p(\lambda)} (3\lambda^2, 2\lambda)$ and
  so $h'(0)=0$.
  
  Conversely, suppose $h \in \mathcal{O}(\mathbb{D})$ satisfies
  $h'(0)=0$.  Then, $f:= h \circ q$ (recall $q:=p^{-1}$) is
  holomorphic on $N\setminus \{(0,0)\}$ because $F(z,w)=h(z/w)$ is
  holomorphic on the set $\{ (z,w) \in \mathbb{D}^2 : |z| < |w| \}$
  which is an open neighborhood of $N\setminus \{(0,0)\}$.  To
  prove $f$ is holomorphic at $(0,0)$, observe first of all that $h$
  can be written as an absolutely convergent power series $h(\lambda)
  = a_0 +a_2 \lambda^2 + a_3 \lambda^3+\cdots$ in some (or any) closed
  disk centered at the origin of radius, say, $r<1$.  Then, for
  $(z,w)$ with $|z|<1, |w|<r^3$,
\[
F(z,w):= a_0 + a_2 w + a_3 z + a_4 w^2 + a_5 z w + a_6 w^3 + \cdots
\]
converges absolutely and extends $f$ (where we are choosing to extend
$(z/w)^k$ to a monomial of the form $zw^m$ or $w^m$---i.e. we want the
power of $w$ to be as large as possible).

Let us emphasize the correspondence just proved:
\begin{observation}
The map given by $f \mapsto f \circ p$ is a bijection from
$\mathcal{O}(N)$ to $\{h\in \mathcal{O}(\mathbb{D}): h'(0)=0\}$ with
inverse given by $h \mapsto h \circ q$, where $p \in
\mathcal{O}(\mathbb{D}, N)$ is $p(\lambda)=(\lambda^3, \lambda^2)$ and
$q=p^{-1}$.
\label{observation1}
\end{observation}

Next, we discuss the complex tangent spaces of $N$.  We can define
$T_x N$ as a subset of $T_x \mathbb{D}^2 \cong \mathbb{C}^2$ in the
following way.  Given $v \in \mathbb{C}^2$, $v \in T_x N$ if and only
if $dG_x v = 0$ for every holomorphic function $G$ defined in a
neighborhood $U$ (in $\mathbb{D}^2$) of $x$ with $G$ identically zero
restricted to $U\cap M$.  Notice that this definition is designed to
make it easy to define the differential of a function $g \in
\mathcal{O}(N)$.

If $x=p(\lambda)=(a,b)\ne (0,0)$ then $T_x N$ is the span of the
vector $(3a,2b)$, because given $G$, as in the previous paragraph, the
function $g:= G \circ p$ is identically zero and so
$0=g'(\lambda)=dG_x(3\lambda^2, 2\lambda)$.  Hence, $dG_x (3a, 2b)=0$.
On the other hand, $h(z,w)=z^2-w^3$ vanishes on $N$ and $dh_x v=0$ if
and only if $v$ is a multiple of $(3a, 2b)$.

At the origin $x=(0,0)$, we have $T_x N =\mathbb{C}^2$, because given
$G$, as above, we have $dG_{(0,0)}=(0,0)$.  This is because the
partial derivatives of $G$ at $(0,0)$ are the coefficients of
$\lambda^3$ and $\lambda^2$ in the identically zero power series for
$G(\lambda^3, \lambda^2)$.  Let us emphasize the above facts:

\begin{observation}
The tangent space $T_{(a,b)} N$ at the point $(a,b) \in N$ with
$(a,b)\ne (0,0)$ can be identified with $\{ \zeta(3a,2b): \zeta \in
\mathbb{C}\} \subset \mathbb{C}^2$.  The tangent space $T_{(0,0)} N$
can be identified with $\mathbb{C}^2$.
\label{observation2}
\end{observation}

\end{subsection}

\end{section}

%%%%%%%%%%%%%%%%%%%%%%%%%%%%%%%%%%%%%%%%%%%%%%%%%%%%%%%%%%%%%%%%%%%%%

\begin{section}{Reduction of theorem \ref{MainThm} to a max problem on $\overline{\mathbb{D}}$}
\label{section:reduction}

As mentioned earlier, we shall compute a formula for $c^{*}_N$ (which
of course gives a formula for $c_N$).

Because of observation \ref{observation1}, we immediately have
\begin{equation}
c^{*}_N(p(\lambda), p(\delta)) = \sup\{m(h(\lambda),h(\delta)): h\in
\mathcal{O}(\mathbb{D},\mathbb{D}), h'(0)=0\}
\label{reduce:eqn}
\end{equation}
As $m$ is invariant under automorphisms of the disk, we may assume
$h(0)=0$ by applying appropriate automorphisms of the disk, since the
condition $h'(0)=0$ is preserved by (post) composition.  Then, by the
Schwarz lemma, $h$ may be written as $h(\zeta)=\zeta^2g(\zeta)$ for
some $g\in \mathcal{O}(\mathbb{D},\overline{\mathbb{D}})$.  At this
stage it is clear that $c^{*}_M(p(0),p(\lambda)) = |\lambda|^2$.  As
$g$ varies over all of
$\mathcal{O}(\mathbb{D},\overline{\mathbb{D}})$, the set of pairs
$(g(\lambda),g(\delta))$ is just the set of all $(A,B)$ satisfying
$m(A,B)\leq m(\lambda, \delta)$.  Hence,
\[
c^{*}_N(p(\lambda),p(\delta)) = \sup\{m(\lambda^2 A, \delta^2 B): A,B
\in \mathbb{D} \text{ with } m(A,B)\leq m(\lambda,\delta)\}
\]
Since $m(\lambda^2 A, \delta^2 B)$ is the modulus of a holomorphic
function in $A$, the above supremum may be taken over all $(A,B)$ with
$m(A,B)=m(\lambda, \delta)$, by the maximum principle.  We may safely
multiply both $A$ and $B$ a unimodular constant and leave $m(\lambda^2
A, \delta^2 B)$ unchanged.  Thus, we can assume there is some
$\alpha\in \mathbb{D}$ such that $A=\phi_\alpha (\lambda)$ and
   $B=\phi_\alpha (\delta)$ (recall $\phi_\alpha$ from \eqref{phi}).

By the preceding discussion we have the following formula for $c^{*}_N$
which gives the desired reduction to a maximization problem.

\begin{prop}
\[
c^{*}_N (p(\lambda),p(\delta)) = \max_{\alpha \in
  \overline{\mathbb{D}}} m(\lambda^2 \phi_{\alpha} (\lambda),\delta^2
  \phi_{\alpha}(\delta))
\]
\label{prop:reduction}
\end{prop}

In particular, the supremum in (\ref{reduce:eqn}) is attained by some
function of the form $h(\zeta) = \zeta^2 \phi_{\alpha}(\zeta)$ where
$\alpha \in \overline{\mathbb{D}}$.  Moreover, if $h$ attains the
supremum in (\ref{reduce:eqn}) and $h(0)=0$, then $h$ is of the same
form (i.e. $h=\zeta^2 \phi_\alpha$) up to multiplication by a
unimodular constant.  As we shall see later, either the supremum will
be obtained with a unique $\alpha \in \mathbb{D}$ or with any $\alpha
\in \partial \mathbb{D}$.  

\end{section}

%%%%%%%%%%%%%%%%%%%%%%%%%%%%%%%%%%%%%%%%%%%%%%%%%%%%%%%%%%%%%%%%%%

\begin{section}{Proof of theorem \ref{MainThm}}
  \label{section:MainThmProof}
To begin, we shall keep $\lambda$ and $\delta$ fixed throughout the
section and define a continuous function, smooth except possibly where
it is zero, $F:\overline{\mathbb{D}} \to [0,1)$ by
\begin{equation}
F(\alpha):= m(\lambda^2 \phi_{\alpha} (\lambda), \delta^2 \phi_{\alpha}
(\delta))
\label{F}
\end{equation}
A couple of things to notice about $F$ are 
\begin{align}
F(\alpha) &< m(\lambda,\delta) \text{ for all } \alpha \in
\overline{\mathbb{D}} \text{ and} \nonumber \\ 
F(\alpha) &=
m(\lambda^2, \delta^2) \text{ for all } \alpha \in \partial \mathbb{D}
\label{Ffact2}.
\end{align}

As in the statement of theorem \ref{MainThm}, we let 
\[
\alpha_0 := \frac{1}{2}\left(\frac{1}{\bar{\lambda}} +
  \lambda+\frac{1}{\bar{\delta}}+\delta \right)
\]

By proposition \ref{prop:reduction}, the following two claims (given
as lemma \ref{lemma1} and lemma \ref{lemma2} below) yield theorem
\ref{MainThm}.  First, $F$ has no local maximum in $\mathbb{D}$ except
possibly $\alpha_0$.  Second, when $|\alpha_0|<1$, $F(\alpha) \leq
F(\alpha_0)$ for all $\alpha$ with $|\alpha|=1$. Before we prove these
facts let us first mention a couple of useful formulas for $F$ whose
proofs we defer to the end of the section.

\begin{claim}\label{claim}

\begin{align}
F(\alpha) &= m(\lambda, \delta)
\left|\frac{(\lambda+\delta)(\alpha+\lambda\delta
    \bar{\alpha}-\lambda-\delta)+\lambda \delta
    (1-|\alpha|^2)}{(1+\lambda \bar{\delta})(1+\lambda \bar{\delta}
    -\bar{\alpha}\lambda-\alpha\bar{\delta})-\lambda
    \bar{\delta}(1-|\alpha|^2)}\right| \label{F1}\\
&= m(\lambda, \delta)
\left|\frac{1 - (\bar{\alpha} - \bar{\alpha}_0 - \bar{\beta_2})(\alpha
    - \alpha_0+\beta_2)}{1 - (\bar{\alpha} - \bar{\alpha}_0 - \bar{\beta_1})(\alpha
    - \alpha_0+\beta_1)}\right| \label{F2}
\end{align}
where
\[
\begin{aligned}
\beta_1 &:= \frac{1}{2}\left(\frac{1}{\bar{\lambda}} -
  \lambda-\frac{1}{\bar{\delta}}+\delta \right) \text{ and}\\
\beta_2 &:= \frac{1}{2}\left(\frac{1}{\bar{\lambda}} -
  \lambda+\frac{1}{\bar{\delta}}-\delta \right)    
\end{aligned}
\]
\end{claim}

%%%% First Lemma %%%%%
\begin{lemma}\label{lemma1}
The function $F$ has no local maximum in $\mathbb{D}$ except possibly
at $\alpha_0$.
\end{lemma}

\begin{proof} Using the formula \eqref{F2}, it suffices to prove the
  function given by
\begin{equation}
G(z)= \left(\frac{F(z+\alpha_0)}{m(\lambda, \delta)}\right)^2 = \left|
  \frac{1-(\bar{z}-\bar{\beta_2})(z+\beta_2)}{1-(\bar{z}-\bar{\beta_1})(z+\beta_1)}
  \right|^2
\label{G}
\end{equation}
has no local max for $|z+\alpha_0|<1$ except possibly at $z=0$.  Some
omitted computations show that $G$ can be written as $G_2/G_1$ where
\begin{equation}
G_k (z)= 1+2|\beta_k|^2-2|z|^2+|z^2-\beta_k^2|^2
\label{Gk}
\end{equation}
for $k=1,2$.

Throughout the following, suppose $z$ is a local maximum satisfying
$0<|z+\alpha_0|<1$.  This implies several things:
\begin{itemize}
\item $0<G(z)<1$, 
\item $z$ is a critical point for $G$, 
\item $\Delta \log G(z) \leq 0$, and
\item $\det {\rm Hess} (\log{G}) \geq 0$ at $z$.  
\end{itemize}
Here ${\rm Hess}$ denotes the matrix of second partial derivatives.
We will prove that all of these conditions cannot be satisfied.

Let us compute all of the derivatives of $G_1$ and $G_2$ up to second
order.  Luckily we can examine $G_1$ and $G_2$ simultaneously.
Writing $z=x+iy$ we have

\[
\begin{aligned}
\partial_z G_k &= -2 \bar{z} +
2z(\bar{z}^2-\bar{\beta_k}^2)\\
\partial_x G_k &= -4x+4{\rm
  Re}[z(\bar{z}^2-\bar{\beta_k}^2)] \\
\partial_y G_k &= -4y-4{\rm
  Im}[z(\bar{z}^2-\bar{\beta_k}^2)] \\
\partial_{xx}^2 G_k &= -4 +4|z|^2+8x^2-4{\rm
  Re}\beta_{k}^2 \\
\partial_{yy}^2 G_k &= -4+4|z|^2+8y^2+4{\rm
  Re}\beta_{k}^2 \\
\partial_{xy}^2 G_k &= 8xy-4{\rm
  Im}\beta_{k}^2
\end{aligned}
\]   

Since $z$ is a critical point for $G$, we have $G_1 \partial_z G_2 -
G_2 \partial_z G_1 = 0$ at $z$.  Neither $G_1$ nor $G_2$ vanish at
$z$, and as a result if $\partial_z G_1 =0$ then $\partial_z G_2 =0$.
But, $\partial_z G_1$ and $\partial_z G_2$ vanish simultaneously only
at 0:
\[
\partial_z G_k = -2 \bar{z} + 2z(\bar{z}^2-\bar{\beta_k}^2)=0
\]
for $k=1,2$ implies $\bar{z}(\beta_1^2-\beta_2^2)=0$, which can only
happen if $z=0$ (because
$\beta_1^2-\beta_{2}^2=-(1-|\lambda|^2)(1-|\delta|^2)/(\bar{\lambda}\bar{\delta})\ne
0$).
Therefore, at $z$
\begin{equation}
\frac{G_2}{G_1} = \frac{\partial_z G_2}{\partial_z G_1}, \qquad
\frac{\partial_x G_1}{G_1} = \frac{\partial_x G_2}{G_2}, \quad {\rm
  and} \quad \frac{\partial_y G_1}{G_1} = \frac{\partial_y G_2}{G_2}
\label{critequations}
\end{equation}
A fact derived from the first of these
equations is
\begin{equation}
\left(\frac{\bar{\beta_1}^2}{G_1} -
  \frac{\bar{\beta_2}^2}{G_2}\right)z^2 =
  |z|^2(1-|z|^2)\left(\frac{1}{G_2}-\frac{1}{G_1}\right)
\label{critequation2}
\end{equation}
and in particular the expression on the left is real.

Using the last two equations in \eqref{critequations}, we can see that
at the critical point $z$ the following equations hold
\[
\begin{aligned}
\partial^2_{xx} \log{G} &= \frac{\partial^2_{xx} G_2}{G_2} -
\frac{\partial^2_{xx} G_1}{G_1} \\ 
&=
(-4+4|z|^2+8x^2)\left(\frac{1}{G_2} - \frac{1}{G_1}\right) + 4{\rm
Re}\left(\frac{\bar{\beta_1}^2}{G_1} -
\frac{\bar{\beta_2}^2}{G_2}\right)\\ 
&= -4[(1-|z|^2)(1-{\rm Re}
(z^2/|z|^2))-2x^2]\left(\frac{1}{G_2}-\frac{1}{G_1}\right)
\end{aligned}
\]
where the last equality follows from (\ref{critequation2}).
Similarly,
\[
\begin{aligned}
\partial^2_{yy} \log{G} &= -4[(1-|z|^2)(1+{\rm Re}(z^2/|z|^2))-2y^2]
\left(\frac{1}{G_2}-\frac{1}{G_1}\right)\\ 
\partial^2_{xy} \log{G} &= 4[2xy+(1-|z|^2){\rm Im}(z^2/|z|^2)]
\left(\frac{1}{G_2}-\frac{1}{G_1}\right)
\end{aligned}
\]

Therefore,
\[
\Delta \log{G} = -8(1-3|z|^2)\left(\frac{1}{G_2}-\frac{1}{G_1}\right)
\]
and as this must be less than or equal to zero at $z$, we see that
$|z|^2 \leq 1/3$.

Finally, we can show that $\det {\rm Hess}(\log{G}) < 0$,
contradicting the fact that $z$ is assumed to be a local maximum.  The
determinant of the Hessian of the logarithm of $G$ (with the positive
factor $16(1/G_2-1/G_1)^2$ omitted) is
\[
\begin{aligned}
  &(1-|z|^2)^2(1-({\rm Re}(z^2/|z|^2))^2)+4x^2y^2
  -2|z|^2(1-|z|^2)\\
  &+2(y^2-x^2)(1-|z|^2){\rm Re}(z^2/|z|^2) \\
  &- 4x^2y^2-4xy(1-|z|^2){\rm Im}(z^2/|z|^2)-(1-|z|^2)^2({\rm
    Im}(z^2/|z|^2))^2
\end{aligned}
\]

Canceling the positive factor $(1-|z|^2)$ and simplifying, we get
$-4|z|^2$ which is indeed negative, as promised.

\end{proof} %% End of proof of first lemma

%%%% Second Lemma %%%%%
\begin{lemma}\label{lemma2} If $|\alpha_0| < 1$, then $F(\alpha) \leq F(\alpha_0)$
  for all $\alpha$ with $|\alpha|=1$.
\end{lemma}
\begin{proof}  Recall  from \eqref{Ffact2} that on the boundary of
  $\overline{\mathbb{D}}$, $F$ is constant and equal to $m(\lambda^2,
  \delta^2)$.  From equation (\ref{F1}) it suffices to prove the
  inequality
\[
\left|\frac{\lambda+\delta}{1+\bar{\lambda}\delta}\right|^2 \leq
\left|\frac{(\lambda+\delta)(\alpha_0+\lambda\delta
    \bar{\alpha}_0-\lambda-\delta)+\lambda \delta
    (1-|\alpha_0|^2)}{(1+\lambda \bar{\delta})(1+\lambda \bar{\delta}
    -\bar{\alpha}_0\lambda-\alpha_0\bar{\delta})-\lambda
    \bar{\delta}(1-|\alpha_0|^2)}\right|^2
\]

Assuming the left hand side above is nonzero (which we can), it
suffices to prove

\begin{align}
&\left|(\alpha_0+\lambda\delta
    \bar{\alpha}_0-\lambda-\delta)+\lambda \delta
    \frac{(1-|\alpha_0|^2)}{\lambda + \delta}\right|^2 \nonumber\\
&- 
\left|(1+\lambda \bar{\delta}
    -\bar{\alpha}_0\lambda-\alpha_0\bar{\delta})-\lambda
    \bar{\delta}\frac{(1-|\alpha_0|^2)}{1+\lambda \bar{\delta}}\right|^2
    \geq 0
\label{inequality}
\end{align}

If we think of the left hand side as
$|A+B|^2-|C+D|^2$=$|A|^2-|C|^2+2{\rm
Re}(A\bar{B}-C\bar{D})+|B|^2-|D|^2$, then first of all $|A|^2-|C|^2$
equals
\[
|\alpha_0+\lambda\delta
    \bar{\alpha}_0-\lambda-\delta|^2-|1+\lambda \bar{\delta}
    -\bar{\alpha}_0\lambda-\alpha_0\bar{\delta}|^2 =
    -(1-|\alpha_0|^2)(1-|\lambda|^2)(1-|\delta|^2)
\]
and using the identities
\[
\begin{aligned}
\alpha_0+\lambda \delta \bar{\alpha}_0-(\lambda+\delta) &=
\frac{\bar{\lambda}+\bar{\delta}}{2\bar{\lambda}\bar{\delta}}
(1+|\lambda \delta|^2)\\
1+\lambda \bar{\delta} -\bar{\alpha}_0\lambda-\alpha_0 \bar{\delta}
&= -\frac{1+\bar{\lambda}\delta}{2\bar{\lambda}\delta}
(|\lambda|^2+|\delta|^2) 
\end{aligned}
\]

we get $2{\rm Re}(A\bar{B}-C\bar{D}) =
(1-|\alpha_0^2)(1-|\lambda|^2)(1-|\delta|^2)$.

Also, using the identity
\begin{equation}
|1+a\bar{b}|^2-|a+b|^2=(1-|a|^2)(1-|b|^2)
\label{identity}
\end{equation}

we see that $|B|^2-|D|^2$ equals
\[
|\lambda
\delta|^2(1-|\alpha_0|^2)^2\frac{(1-|\lambda|^2)(1-|\delta|^2)}{|\lambda+\delta|^2|1+\lambda\bar{\delta}|^2}
\]

Summing this all up, we see that proving (\ref{inequality}) amounts to
showing
\[
|\lambda
 \delta|^2(1-|\alpha_0|^2)^2\frac{(1-|\lambda|^2)(1-|\delta|^2)}{|\lambda+\delta|^2|1+\lambda\bar{\delta}|^2}
 \geq 0
\]
which is certainly true.

\end{proof}

This concludes the proof of theorem \ref{MainThm}.  As promised, a
slightly nicer formula for $c^{*}_N (p(\lambda), p(\delta))$ is
\begin{prop} If $\lambda, \delta \in \mathbb{D}$ are nonzero, then
\[
c^{*}_N (p(\lambda), p(\delta)) = \begin{cases} m(\lambda^2, \delta^2)
  & \text{ if } |\alpha_0| \geq 1 \\
m(\lambda, \delta) \frac{1+|\beta_2|^2}{1+|\beta_1|^2} & \text{ if }
  |\alpha_0| < 1
\end{cases}
\]
\label{prop:nicer}
\end{prop}

This follows from proposition \ref{prop:reduction}, the definition of
$F$ (namely equation \eqref{F}), formula \eqref{F2} for $F$, and the
above lemmas.

We conclude this section with the proof of claim \ref{claim}.

\begin{proof}[Proof of Claim \ref{claim}:]
We start from equation (\ref{F}).  Observe that
\begin{align}
F(\alpha) &=
\left|\frac{\lambda^2\frac{\alpha-\lambda}{1-\bar{\alpha}\lambda} -
    \delta^2\frac{\alpha-\delta}{1-\bar{\alpha}\delta}}
  {1-\lambda^2\bar{\delta}^2\frac{\alpha-\lambda}{1-\bar{\alpha}
        \lambda}\frac{\bar{\alpha}-\bar{\delta}}{1-\alpha\bar{\delta}}}\right| \nonumber\\
&= \left|\frac{\lambda^2(\alpha-\lambda)(1-\bar{\alpha}\delta) -
    \delta^2(\alpha-\delta)(1-\bar{\alpha}\lambda)}
  {(1-\bar{\alpha}\lambda) (1-\alpha \bar{\delta}) - \lambda^2
    \bar{\delta}^2
    (\alpha-\lambda)(\bar{\alpha}-\bar{\delta})}\right|\nonumber\\
&= \left|\frac{\alpha(\lambda^2-\delta^2)
    -(\lambda^3-\delta^3)-|\alpha|^2\lambda \delta
    (\lambda-\delta)+\lambda \delta (\lambda^2-\delta^2)\bar{\alpha}}
  {1-\lambda^3 \bar{\delta}^3-\bar{\alpha} \lambda (1-\lambda^2
    \bar{\delta}^2)-\alpha \bar{\delta} (1-\lambda^2
    \bar{\delta}^2)+|\alpha|^2\lambda \bar{\delta} (1-\lambda
    \bar{\delta})}\right|\nonumber\\
&= m(\lambda, \delta)
\left|\frac{\alpha(\lambda+\delta)-(\lambda^2+\lambda \delta +
    \delta^2)-|\alpha|^2 \lambda \delta + \lambda \delta (\lambda
    +\delta)\bar{\alpha}} {1+\lambda
    \bar{\delta}+\lambda^2\bar{\delta}^2 -\bar{\alpha}\lambda
    (1+\lambda \bar{\delta})-\alpha \bar{\delta} (1+\lambda
    \bar{\delta})+|\alpha|^2 \lambda \bar{\delta}}\right| \label{lasteqn}\\
&= m(\lambda, \delta)\left|\frac{\alpha(\lambda+\delta)+\lambda
    \delta (\lambda + \delta) \bar{\alpha}- (\lambda +
    \delta)^2+\lambda \delta (1-|\alpha|^2)} {(1+\lambda
    \bar{\delta})^2 - \bar{\alpha} \lambda (1+\lambda \bar{\delta})-
    \alpha \bar{\delta} (1+\lambda \bar{\delta}) - (1-|\alpha|^2)
    \lambda \bar{\delta}}\right| \nonumber
\end{align}
and from here it is easy to get (\ref{F1}).

Secondly, to prove (\ref{F2}), we start from (\ref{lasteqn}):
\[
\begin{aligned}
F(\alpha) &= m(\lambda, \delta)\left| \frac{\lambda \delta - (\lambda
    \delta \bar{\alpha} - (\lambda + \delta)) (\alpha - (\lambda +
    \delta))} {\lambda \bar{\delta} - (\bar{\alpha} \lambda
    -(1+\lambda \bar{\delta})) (\alpha \bar{\delta} - (1+\lambda
    \bar{\delta}))}\right|\\
&= m(\lambda, \delta) \left| \frac{1-(\bar{\alpha} - (1/\delta +
    1/\lambda)) (\alpha- (\lambda + \delta))} {1- (\bar{\alpha} -
    (1/\lambda +\bar{\delta})) (\alpha - (1/\bar{\delta} + \lambda))}
    \right|
\end{aligned}
\]
and this equals (\ref{F2}) because of the identities:
\[
\begin{aligned}
\bar{\alpha}_0 + \bar{\beta_2} &=& \frac{1}{\lambda} +
\frac{1}{\delta} \\
\alpha_0-\beta_2 &=& \lambda + \delta \\
\bar{\alpha}_0 + \bar{\beta_1} &=& \frac{1}{\lambda} +\bar{\delta}\\
\alpha_0 - \beta_1 &=& \lambda + \frac{1}{\bar{\delta}}
\end{aligned}
\]

\end{proof} % end of proof of claim

\end{section} % end of proof of main theorem

%%%%%%%%%%%%%%%%%%%%%%%%%%%%%%%%%%%%%%%%%%%%%%%%%%%%%%%%%%%%%

\begin{section}{The infinitesimal Carath\'eodory pseudodistance}
\label{section:infproof}
In this section we prove theorem \ref{thm2}, our formula for the
Carath\'eodory pseudometric.

The Carath\'eodory pseudometric at the origin and a vector
$v=(v_1,v_2)\in \mathbb{C}^2$ is
\[
C_N((0,0); v) = \sup\{|df_{(0,0)} v| : f\in \mathcal{O}(N,\mathbb{D})
\text{ and } f(0,0)=0\}
\]
Any $f$ as above satisfies $f(\lambda^3, \lambda^2) = \lambda^2
g(\lambda)$ for some $g \in \mathcal{O}(\mathbb{D},
\overline{\mathbb{D}})$ (see the beginning of section
\ref{section:reduction}).  Also, the partial derivative of $f$ with
respect to the first variable at the origin is just $g'(0)$ and the
partial derivative of $f$ with respect to the second variable at the
origin is $g(0)$ (see section \ref{section:discussion}).  Therefore,
\[
C_N((0,0); v) = \sup\{ |v_1 g'(0)+v_2 g(0)| : g \in
\mathcal{O}(\mathbb{D},\overline{\mathbb{D}})\}
\]
The set of pairs $(g'(0), g(0))$ as $g$ varies over
$\mathcal{O}(\mathbb{D},\overline{\mathbb{D}})$ is really just the
pairs $(A,B)$ where $|A|+|B|^2 \leq 1$, by the Schwarz-Pick Lemma.
With suitable choices for the arguments of $A$ and $B$, we can reduce
the problem to maximizing $|v_1|s+|v_2|t$ over all $s,t \in [0,1]$
satisfying $s+t^2 \leq 1$.  The function we are maximizing is linear,
so the maximum occurs on the boundary.  Therefore, the problem is just
a matter of finding the maximum of $|v_1|(1-t^2)+|v_2|t$ for $0\leq t
\leq 1$.  So, by calculus,
\[
C_N ((0,0); v) = \begin{cases} |v_2| & \text{ if } |v_2| \geq 2 |v_1| \\
                               \frac{4|v_1|^2+|v_2|^2}{4|v_1|} &
                               \text{ if } |v_2| < 2 |v_1| \end{cases}
\]
as desired.

Next, let $x=(a,b) \in N\setminus \{(0,0)\}$ and define $v = (3a,2b)$.
The Carath\'eodory pseudometric at $(a,b)$ is 
\[
C_N(x; v) = \sup\left\{ \frac{|df_x v|}{1-|f(x)|^2} : f \in
  \mathcal{O}(N, \mathbb{D})\right\}
\]
If we set $\lambda = a/b$ and $h=f \circ p$, then $v = \lambda
(3\lambda^2,2\lambda)$ and since $df_x (3\lambda^2, 2\lambda) =
h'(\lambda)$ we see that
\[
C_N(x;v) = |\lambda| \sup\{\rho(h(\lambda); h'(\lambda)) : h \in
\mathcal{O}(\mathbb{D},\mathbb{D}) \text{ and } h'(0)=0\}
\]
By (post) composing $h$ with an automorphism of the unit disk (which
is allowed by invariance properties of $\rho$), we can assume $h(0)=0$
and therefore $h$ has the form $h(\zeta)=\zeta^2 g(\zeta)$ for some $g
\in \mathcal{O}(\mathbb{D},\overline{\mathbb{D}})$.  Hence,
\[
C_N(x;v) = |\lambda| \sup\left\{ \frac{|\lambda^2 g'(\lambda)+2\lambda
    g(\lambda)|}{1-|\lambda|^4 |g(\lambda)|^2} : g \in
    \mathcal{O}(\mathbb{D},\overline{\mathbb{D}}) \right\}
\]
Like before, $(g'(\lambda), g(\lambda))$ varies over all pairs $(A,B)$
satisfying $|A|(1-|\lambda|^2) \leq 1-|B|^2$.  This reduces the
problem to maximizing
\[
\frac{|\lambda|^2 s+2|\lambda|t}{1-|\lambda|^4 t^2}
\]
over the set of non-negative $s,t$ satisfying $t^2+s(1-|\lambda|^2)
\leq 1$.  It is easy to check that the maximum always occurs when
$t=1$ and $s=0$.  Since $\lambda^2=b$ we see that
\[
C_N(x;v) = \frac{2|b|}{1-|b|^2}.
\]
\end{section}

\begin{section}{Proof of theorem \ref{prop:interp}}
\label{section:proveinterp}
By precomposing all functions with $\phi_{z_3}$ we may assume $z_3=0$
in theorem \ref{prop:interp}.  Then, all functions of interest will
correspond to functions in $\mathcal{O}(N,\mathbb{D})$, and therefore
it is clear that if there is a function $h \in
\mathcal{O}(\mathbb{D},\mathbb{D})$ which satisfies the $h'(0)=0,
h(z_i)=w_i$ for $i=1,2$, then the inequality (\ref{mixedcp-ineq})
holds (by theorem \ref{MainThm} and the definition of Carath\'eodory
pseudodistance).

On the other hand, if the inequality (\ref{mixedcp-ineq}) holds (again
with $z_3=0$), then pick a function $f \in \mathcal{O}(N,\mathbb{D})$
with
\[
\rho (f(p(z_1)), f(p(z_2))) = c_N (p(z_1), p(z_2))
\]
(we know such a function exists by the formula for $c_N$) and then set
$h:=f \circ p \in \mathcal{O}(\mathbb{D},\mathbb{D})$.  The function
$h$ satisfies $\rho (w_1, w_2) \leq \rho (h(z_1), h(z_2))$ and by
composing $h$ with an appropriate function we can find a function $F
\in \mathcal{O}(\mathbb{D},\mathbb{D})$ with $F(z_1)=w_1$,
$F(z_2)=w_2$, and $F'(0)=0$.

To prove the last part of theorem \ref{prop:interp}, suppose $F$
satisfies the interpolation problem and equality in
(\ref{mixedcp-ineq}).  Then, $h:=\phi_{F(0)} \circ F$ satisfies
equality as well.  Hence, when
\[
\alpha_0 :=
\frac{1}{2}\left(\frac{1}{\bar{z}_1}+z_1+\frac{1}{\bar{z}_2} + z_2
\right)
\]
is in the disk, $h(\lambda)$ is of the form $\mu \lambda^2
\phi_{\alpha_0}(\lambda)$ where $\mu$ is a unimodular constant, and when
$\alpha_0 \notin \mathbb{D}$, $h(\lambda)$ is of the form $\mu \lambda^2$
(again with $\mu \in \partial\mathbb{D}$).  But, $\mu$ and $F(0)$ are
uniquely determined by the fact that $w_i = \phi_{F(0)}( h(z_i))$ for
$i=1,2$ since $h(z_1)$ and $h(z_2)$ must be distinct.  So, there
exists a unique automorphism of the disk $\psi$ such that
\[
F(\lambda) = \begin{cases} \psi( \lambda^2 \phi_{\alpha_0}(\lambda)) & \text{
    if } \alpha_0 \in \mathbb{D} \\
                         \psi( \lambda^2) & \text{ if } \alpha_0 \notin
    \mathbb{D} \end{cases}
\]
In the first case, $F$ is a Blaschke product of order three and in the
second a Blaschke product of order two.
\end{section}

%%%%%%%%%%%%%%%%%%%%%%%%%%%%%%%%%%%%%%%%%%%%%%%%%%%%%%%%%%%%%%%%%%%

\begin{section}{Proof of Extension Theorem}
\label{section:proveextend}
In this section we prove Theorem \ref{thm:extend}.  

First, we need to define a few basic notions. Let $X$ be a set.  A
self-adjoint function $F: X \times X \to \mathbb{C}$ (i.e.
$F(x,y)=\overline{F(y,x)}$) is \emph{positive semi-definite} if for
every positive integer $n$ and every finite subset $\{x_1,x_2,\dots,
x_n\} \subset X$ the $n\times n$ matrix with entries $F(x_i, x_j)$ is
positive semi-definite.  For example, by the Pick interpolation
theorem the function $F: \mathbb{D} \times \mathbb{D} \to \mathbb{C}$
given by
\[
F(\lambda, \delta) = \frac{1-g(\lambda)\overline{g(\delta)}}{1-\lambda
  \bar{\delta}}
\]
is positive semi-definite for any $g \in
\mathcal{O}(\mathbb{D},\overline{\mathbb{D}})$.  

The Pick interpolation theorem on the bidisk (see \cite{AM} page 180)
can be stated as a theorem about extensions of bounded analytic
functions in the following way.  Given a subset $X$ of the bidisk, and
a function $\psi: X \to \mathbb{D}$ there exists $\Psi \in
\mathcal{O}(\mathbb{D}^2, \mathbb{D})$ with $\Psi\mid_X = \psi$
if and only if there exist positive semi-definite functions $\Delta$
and $\Gamma$ on $X\times X$ such that for each $z=(z_1, z_2),
w=(w_1,w_2) \in X$
\[
1-\psi(z)\overline{\psi(w)} = \Gamma (z,w) (1-z_1 \bar{w}_1) + \Delta (z,w)
(1-z_2 \bar{w}_2)
\]
We should mention that the portion of this theorem which we shall use
(namely sufficiency) has a quite simple proof---it is an application
of the so-called ``lurking isometry'' technique.

To prove theorem \ref{thm:extend} suppose $f \in
\mathcal{O}(N,\mathbb{D})$ and $f(0,0)=0$.  Then, as in earlier
arguments, $(f\circ p)(\lambda)=f(\lambda^3, \lambda^2) = \lambda^2 g(\lambda)$
for some $g \in \mathcal{O}(\mathbb{D}, \overline{\mathbb{D}})$.  For
any $\delta, \lambda \in \mathbb{D}$, we have
\[
\begin{aligned}
2-&f(p(\lambda))\overline{f(p(\delta))} = (1-\lambda^3 \bar{\delta}^3)\\
& + \left( 1 + \lambda^2
  \bar{\delta}^2\frac{1-g(\lambda)\overline{ g(\delta)}}{1-\lambda
  \bar{\delta}} + \frac{\lambda^3 \bar{\delta}^3 g(\lambda)\overline{
  g(\delta)}}{ 1-\lambda^2 \bar{\delta}^2} \right) (1-\lambda^2
  \bar{\delta}^2)
\end{aligned}
\]

Therefore, for $z=(z_1,z_2), w=(w_1,w_2) \in N$ 
\begin{equation}
2-f(z)\overline{f(w)} = \Gamma (z,w) (1-z_1 \bar{w}_1) + \Delta (z,w)
(1-z_2\bar{w}_2) 
\label{extendf}
\end{equation}

where $\Gamma (z,w) = 1$ and
\[
\Delta (z,w) =1+z_1 \bar{w}_1 \frac{1-
  g(q(z))\overline{g(q(w))}}{1-q(z) \overline{q(w)}} + \frac{z_2
  \bar{w}_2 g(q(z))\overline{g(q(w))}}{ 1-z_1 \bar{w}_1}
\]
(recall $q(z)=z_1/z_2$ for $z\ne (0,0)$ and $q(0,0)=0$).  Now,
$\Gamma$ is clearly positive semi-definite, and $\Delta$ is positive
semi-definite because of the fact that positive semi-definite
functions are closed under addition and multiplication (by the Schur
product theorem) and by the Pick interpolation theorem on the disk
(applied to $g$).  This proves $f$ has an extension to the bidisk with
supremum norm at most $\sqrt{2}$ (by dividing through (\ref{extendf})
by 2).

To prove any holomorphic function $f \in \mathcal{O}(N,\mathbb{D})$
(regardless of its value at the origin) can be extended to the bidisk
with supremum norm at most $2\sqrt{2}+1$, simply apply the result just
proved to $(f-f(0))/2$.

Finally, the function
\[
h(\lambda) = \lambda^2 \frac{0.5-\lambda}{1-0.5\lambda}
\]
corresponds to a function $f \in \mathcal{O}(N,\mathbb{D})$ with
$f(\lambda^3, \lambda^2) = h (\lambda)$.  The partial derivatives of
$f$ at $(0,0)$ are just the coefficients of $\lambda^3$ and
$\lambda^2$ in the power series for $h$; i.e. they are $-0.75$ and
$0.5$.  Suppose $F$ is a bounded holomorphic extension of $f$ to the
bidisk with sup norm $R$.  Then, by the Schwarz lemma on the bidisk
\[
0.75/R +0.5/ R \leq 1
\]
which implies $R \geq 5/4$, as desired.

\subsection*{Acknowledgements}
Thanks to Marek Jarnicki and Peter Pflug for posing the problem which
led to this paper, to my advisor, John M{\raise .45ex \hbox{c}}Carthy,
for numerous suggestions, and to Nikolai Nikolov for pointing out an
error in an earlier version of this paper.

\end{section}

%%%%%%%%%%%%%%%%%%%%%%%%%%%%%%%%%%%%%%%%%%%%%%%%%%%%%%%%%%%%%%%%%%


\begin{thebibliography}{1}
\bibitem{AM} J.\! Agler and J.E.\! M{\raise .45ex \hbox{c}}Carthy, \emph{Pick
    Interpolation and Hilbert Function Spaces}, A.M.S., Providence,
    RI, 2002.
    
  \bibitem{AM-03} J.\! Agler and J.E.\! M{\raise .45ex
      \hbox{c}}Carthy, \emph{Norm Preserving Extensions of Holomorphic
      Functions from Subvarieties of the Bidisk}, Ann. of Math.,
      \textbf{157} (2003), 289-312.
    
  \bibitem{CF} C.\! Carath\'eodory and L.\! Fej\'er, \emph{\"Uber den
      Zusammenhang der Extremen von harmonischen Funktionen mit ihren
      Koeffizienten und \"uber den Picard-Landauschen Satz}, Rend.
    Circ. Mat. Palermo. \textbf{32} (1911), 218-239.

  \bibitem{FF} C.\! Foia\c{s} and A.E.\! Frazho, \emph{The commutant
  lifting approach to interpolation problems}, Birkh\"auser, Basel,
  1990.
    
  \bibitem{GUN} R.\! Gunning, \emph{Introduction to Holomorphic
      Functions of Several Variables, Vol II}, Wadsworth Inc.,
    California, 1990.
    
  \bibitem{GUN3} R.\! Gunning, \emph{Introduction to Holomorphic
      Functions of Several Variables, Vol III}, Wadsworth Inc.,
      California, 1990.
    
  \bibitem{IK} A.\! Isaev and S.\! Krantz, \emph{Invariant Distances
      and Metrics in Complex Analysis}, Notices of the American Math.
      Society, \textbf{47} (2000), 546-553.

\bibitem{J-P93} M.\! Jarnicki and P.\! Pflug, \emph{Invariant
    Distances and Metrics in Complex Analysis}, de Gruyter Expositions
    in Mathematics 9, Walter de Gruyter, 1993.
    
  \bibitem{J-P} M.\! Jarnicki and P.\! Pflug, \emph{Invariant
      Distances and Metrics in Complex Analysis - Revisited},
      Diss. Math.  \textbf{430} (2005), 1-192.
    
  \bibitem{KOB} S.\! Kobayashi, \emph{Hyperbolic Complex Spaces},
    Springer Verlag, 1998.
    
  \bibitem{LEM} L.\! Lempert, \emph{La m\'etrique de Kobayashi et la
      repr\'esentation des domaines sur la boule}, Bull. Soc. Math.
    France \textbf{109} (1981), 427-474.

\bibitem{NEV} R.\! Nevanlinna, \emph{\"Uber beschr\"ankte Funktionen},
  Ann. Acad. Sci. Fenn. Ser. A \textbf{13} (1919), no. 1.
    
\bibitem{PICK} G.\! Pick, \emph{\"Uber die Beschr\"ankungen analytischer
    Funktionen, welche durch vorgegebene Funktionswerte bewirkt
    werden}, Math. Ann.  \textbf{77} (1916), 7-23.
    
\bibitem{P-H} P.L.\! Polyakov and G.M.\! Khenkin, \emph{Integral Formulas
    for the $\bar{\partial}$-Equation and an Interpolation Problem in
    Analytic Polyhedra}, Trans. Moscow Math. Soc.  \textbf{53} (1991),
  135-175.

  \bibitem{WAL} J.\! Wallis, \emph{Arithmetica Infinitorum}, 1655.

\end{thebibliography}
\end{document}